\documentclass[a4paper,10pt]{amsart}

\usepackage{amsmath,amsthm,amstext,amsfonts,float,amssymb,mathrsfs}
\usepackage{mathptmx}
\usepackage[all]{xy}

\title[A symplectic construction of Calabi's extremal K\"ahler metrics]
{A symplectic construction of Calabi's extremal K\"ahler metrics on
the blow-up of $\pmb{\bbC\bbP^n}$ at one point}

\author{Aleksis Raza}

\address{Department of Mathematics, 180 Queens Gate, Imperial College London SW7 2BT UK}

\email{aleksis.raza@imperial.ac.uk}

\newtheorem{defin}{Definition}

\newtheorem{theorem}[defin]{Theorem}

\newtheorem{coro}[defin]{Corollary}

\newtheorem{example}[defin]{Example}

\newcommand{\bbR}{\mathbb{R}}
\newcommand{\bbC}{\mathbb{C}}

\newcommand{\bbP}{\mathbb{P}}

\newcommand{\tts}{{\widetilde{s}}}

\newcommand{\pon}{\bbR^n_{\geq 0}}

\newcommand{\U}{{\rm U}}
\newcommand{\SU}{{\rm SU}}

\renewcommand{\ln}{\log}
\renewcommand{\geq}{\geqslant}

\begin{document}

\maketitle

\begin{abstract}
We apply a local differential geometric framework from K\"ahler
toric geometry to (re)construct Calabi's extremal K\"ahler metrics
on $\bbC\bbP^n$ blown-up at a point from data on the moment
polytope. This note is an addendum to
\cite{raz04:app_guil_abr_non_abel_group_act}.
\end{abstract}

\section{Recollection of relevant results}

Using a construction from symplectic toric geometry (see \S 1 of
\cite{raz04:app_guil_abr_non_abel_group_act}) we proved the
following
\begin{theorem}[\cite{raz04:scal_cur_mul_free_act}]
Suppose $T^n$ acts on $\bbC^n$ (with its fixed standard symplectic
structure) via the standard linear action so that the moment
polytope $\Delta_{\bbC^n}$ of the moment map corresponding to this
action is the positive orthant $\pon$ in $\bbR^n\cong(\bbR^n)^*$
with standard symplectic (action) coordinates
$(x_1,\dots,x_n)=\frac{1}{2}(|z_1|^2,\dots,|z_n|^2)$, where
$(z_1,\dots,z_n)$ are the standard complex coordinates on $\bbC^n
\setminus \{0\}$. Now consider the $\SU(n)$-invariant K\"ahler
metric $h_f$ on $\bbC^n \setminus \{0\}$ determined by the standard
complex structure on $\bbC^n$ and the K\"ahler potential $f(s)$
where $s=\sum_{i=1}^n|z_i|^2$. Then there exists a K\"ahler isometry
between the K\"ahler manifold $(\bbC^n \setminus \{0\},h_f)$ and the
K\"ahler manifold $(\bbC^n \setminus \{0\},h_g)$ where $h_g$ is the
$\SU(n)$-invariant K\"ahler metric on $\bbC^n \setminus \{0\}$
determined by the standard symplectic structure on $\bbC^n$ and the
symplectic potential
\begin{equation}\label{spp}
g(x)=\frac{1}{2}\left(\sum_{i=1}^n x_i \ln x_i + F(t)\right)
\end{equation}
which is a smooth function defined on the interior
$\Delta^\circ_{\bbC^n}=\bbR^n_{>0}$ of $\Delta_{\bbC^n}$, where
\[
t=\sum_{i=1}^nx_i=2sf'(s)
\]
and
\[
F(t)= t \ln \left(s(t)t^{-1}\right) - 2f(s(t))
\]
which we call the $t$-potential of \eqref{spp}. Furthermore, the
scalar curvature of $h_g$ (and hence of $h_f$) is given (in the
coordinates $x$) by
\begin{equation}\label{ss}
S(g)=t^{1-n}\left(t^{n+1}F''(t)(1+tF''(t))^{-1}\right)''.
\end{equation}
Conversely, any function of the form $g(x)$ on
$\Delta^\circ_{\bbC^n}$ with $F(t)$ satisfying $ F''(t)>-t^{-1}$
determines such a K\"ahler metric.
\end{theorem}
Note that $\SU(n)$-invariant K\"ahler metrics on $\bbC^n \setminus
\{0\}$ are in fact $\U(n)$ (hence $T^n$)-invariant which is why the
above result comes about using ideas from toric geometry. Motivated
by the discussion in \S 3 of \cite{abr03:kah_geom_tor_man_sym_coor}
and the example in \S 6 of \cite{abr98:kah_geom_tor_var_ext_met} we
shall now use this theorem to (re)construct a family of extremal
K\"ahler metrics on certain $\bbC\bbP^1$ bundles over
$\bbC\bbP^{n-1}$ originally identified by Calabi
\cite{cal82:ext_kah_met}.

\section{Calabi's extremal K\"ahler metrics on
$\pmb{\widehat{\bbC\bbP}^n}$}

According to Abreu
\begin{theorem}[\cite{abr98:kah_geom_tor_var_ext_met}]
A toric K\"ahler metric determined by a symplectic potential $g(x)$
is extremal if and only if
\begin{equation}\label{excon}
\frac{\partial S}{\partial x_i}={\rm constant}
\end{equation}
for $i=1,\dots,n$ i.e. its scalar curvature
\begin{equation}\label{scalt}
S(g) = - \frac{1}{2}\sum^{n}_{i,j=1} \frac{\partial^2
G^{ij}}{\partial x_i \partial x_j}
\end{equation}
is an affine function of $x$. Here $G^{ij}$ is the $(i,j)$th entry
of inverted hessian matrix $G^{-1}$ of $g$.
\end{theorem}

In \S 3, pp.278-88 of \cite{cal82:ext_kah_met} Calabi constructed
extremal K\"ahler metrics of non-constant scalar curvature on
certain $\bbC\bbP^1$ bundles over $\bbC\bbP^{n-1}$. Since the total
space of these bundles is $\bbP(\mathscr{O}(-1)\oplus \bbC)$ i.e.
the projectivization of the direct sum of the tautologous line
bundle over $\bbC\bbP^{n-1}$ with the trivial line bundle, one may
equally regard these manifolds as the blow-up of $\bbC\bbP^n$ at a
point, $\widehat{\bbC\bbP}^n$. Calabi constructed these metrics
using a scalar curvature formula ((3.9) in \S 3 of
\cite{cal82:ext_kah_met}) for K\"ahler metrics on $\bbC^n\setminus
\{0\}$ (considered as an open set in these bundles) invariant under
$\U(n)$ (the maximal compact subgroup of the group of complex
automorphisms of these bundles) and then imposing the appropriate
boundary conditions. This scalar curvature formula was
\begin{equation}\label{cal}
S(f)= (n-1)v'(\tts) \left(f'(\tts)\right)^{-1} +v''(\tts)
\left(f''(\tts)\right)^{-1}
\end{equation}
where $ v(\tts)= n\tts-(n-1)\ln f'(\tts) - \ln f''(\tts)$ and $
\tts=\ln s$. Using the symplectic (action) coordinate
$t=\sum_{i=1}^nx_i$ \eqref{cal} becomes of the extremely simple form
$\eqref{ss}$. We will now (re)construct Calabi's metrics using
\eqref{ss} and working on the Delzant moment polytope
$\Delta_{\widehat{\bbC\bbP}^n}$ of the symplectic toric manifold
$\widehat{\bbC\bbP}^n$. Our motivation for doing this came from the
discussion in \S 3 of \cite{abr03:kah_geom_tor_man_sym_coor}.

\subsection{The construction}

$\Delta_{\widehat{\bbC\bbP}^n}$ consists of $n+2$ facets and these
are determined by the affine functions
\begin{equation}\label{li}
l_i(x)=\left\{\begin{array}{cl}x_i & i=1,\dots,n
\\ t-a & i=n+1 \\ b-t & i=n+2 \end{array}\right.
\end{equation}
where $0<a<b$. The constants $a,b$ determine the cohomology class of
the extremal K\"ahler metrics we are about to construct. The
extremal condition \eqref{excon} in the $\U(n)$-invariant scenario
becomes $dS/dt={\rm constant}$ i.e. $S=At+B$ for constants $A,B$.
Thus we have \eqref{ss}$=At +B$. Solving this for $F''$ gives
\begin{equation}\label{text}
F''(t)=\frac{pt^{n-1}}{pt^n-\alpha} - \frac{1}{t}
\end{equation}
where $p=n(n+1)(n+2)$ and
\begin{equation}\label{alpha}
\alpha(t)=nAt^{n+2}+(n+2)Bt^{n+1}+p(Ct+D)
\end{equation}
for constants $C,D$. \eqref{text} is a general formula for the
$t$-potential of any extremal metric on $\bbC^n\setminus\{0\}$. We
denote by $F _{\widehat{\bbC\bbP}^n}$ the $t$-potential of the
extremal K\"ahler metrics we are after. By \eqref{spp} the
symplectic potential of these metrics is
$2g_{\widehat{\bbC\bbP}^n}=\sum_{i=1}^n x_i \ln x_i +
F_{\widehat{\bbC\bbP}^n}(t)$ which is
\begin{equation}\label{symp}
g_{\widehat{\bbC\bbP}^n}=\frac{1}{2}\left[\sum_{i=1}^{n+2} l_i \ln
l_i + h_{\widehat{\bbC\bbP}^n}(t)\right]
\end{equation}
in the form (2.9) of Abreu's Theorem 2.8 in
\cite{abr03:kah_geom_tor_man_sym_coor} with the $l_i$ given by
\eqref{li}. Therefore $h_{\widehat{\bbC\bbP}^n}(t) =
F_{\widehat{\bbC\bbP}^n}(t) - \sum_{i=n+1}^{n+2} l_i \ln l_i.$
Differentiating this twice gives
\begin{equation}\label{symp2}
h''_{\widehat{\bbC\bbP}^n}(t) = F''_{\widehat{\bbC\bbP}^n}(t) -
\frac{b-a}{(t-a)(b-t)}=\frac{pt^{n+1}-2apt^n+abpt^{n-1} -c \alpha}
{(pt^n-\alpha)(t-a)(t-b)} -\frac{1}{t}
\end{equation}
where $c=b-a$. Thus $a<t<b$ and the boundary conditions are that
as $ t \downarrow a$ and $t \uparrow b$ then
\begin{equation}\label{approx}
pt^{n+1}-2apt^n+abpt^{n-1} -c \alpha \approx
(pt^n-\alpha)(t-a)(b-t).
\end{equation}
Consider the former condition. Let $t=a+\epsilon$ for small
$\epsilon >0$ and ignore terms of $O(\epsilon^2)$ and above. Then
\eqref{alpha} gives $\alpha=\alpha_a + O(\epsilon^2)$ where
$\alpha_a=X_a+\epsilon Y_a$ with
\begin{equation}\label{x1}
X_a=na^{n+2}A+(n+2)a^{n+1}B+paC+pD
\end{equation}
and
\begin{equation}\label{y1}
Y_a=n(n+2)Aa^{n+1}+ (n+1)(n+2)Ba^n +pC.
\end{equation}
Now $(t-a)(b-t) = -c\epsilon +O(\epsilon^2)$ and $pt^n-\alpha_a =
pa^n+pna^{n-1}\epsilon -\alpha_a+O(\epsilon^2)$. Therefore $
(pt^n-\alpha)(t-a)(b-t) = -cpa^n\epsilon +c X_a\epsilon +
O(\epsilon^2)$. An analogous calculation shows that $
pt^{n+1}-2apt^n+abpt^{n-1} = cpa^n +  cpa^{n-1}(n-1)\epsilon +
O(\epsilon^2)$. Hence \eqref{approx} gives $ cpa^n +
cpa^{n-1}(n-1)\epsilon -cX_a-cY_a\epsilon = -cpa^n\epsilon +c
X_a\epsilon$ which rearranges into $ pa^n + (pa^{n-1}(n-1) +
pa^n)\epsilon = X_a + (X_a+Y_a)\epsilon$. Comparing coefficients
shows that
\begin{equation}\label{x12}
X_a= pa^n
\end{equation}
and $ X_a+Y_a = pa^{n-1}(n-1) + pa^n $ i.e.
\begin{equation}\label{y12}
Y_a= (n-1)pa^{n-1}.
\end{equation}
The calculations for the second condition $t \uparrow b$ follow in
similar way i.e. we consider $t=b-\epsilon$ for small $\epsilon>0$.
By \eqref{alpha} we have that $\alpha=\alpha_b+O(\epsilon^2)$ with
$\alpha_b=X_b-\epsilon Y_b$ where $X_b$ and $Y_b$ are \eqref{x1} and
\eqref{y1}, respectively, with $a$ replaced by $b$. Also
$(t-a)(b-t)=c\epsilon +O(\epsilon^2)$ so that $
(pt^n-\alpha)(t-a)(b-t) = cpb^n\epsilon - cX_b\epsilon$.
Furthermore, $ pt^{n+1}-2apt^n+abpt^{n-1} =
pcb^n-(n+1)pcb^{n-1}\epsilon + O(\epsilon^2)$. By \eqref{approx} we
get $pcb^n-(n+1)pcb^{n-1} - cX_b+cY_b\epsilon = cpb^n\epsilon -
cX_b\epsilon $ which rearranges to give $ pcb^n-((n+1)pcb^{n-1} +
cpb^n)\epsilon = cX_b - c(X_b +Y_b)\epsilon$. Comparing coefficients
gives
\begin{equation}\label{x22}
X_b= pb^n
\end{equation}
and $ X_b+Y_b = ((n+1)pb^{n-1} + pb^n) $ i.e.
\begin{equation}\label{y22}
Y_b=(n+1)pb^{n-1}.
\end{equation}
Thus \eqref{x12}-\eqref{y22} give us four linear equations in the
unknowns $A,B,C,D$. Solving these simultaneously reveals
\begin{equation}\label{soln}
\begin{array}{c}
A=\frac{(n+1)(n+2)((ab)^{n-1}(na^2(n+1) +nb^2(n-1)
-2ab(n^2-1))-2a^{2n})}
{(ab)^n(2n(n+2)ab-(a^2+b^2)(n+1)^2)+a^{2(n+1)} + b^{2(n+1)}}, \\
\\
B=\frac{n(n+1)((ab)^{n-1}(a^2(nb(n+2)-a(n+1)^2) +
b^2(b(1-n^2)+a(n^2-4)))+3a^{2n+1}+b^{2n+1})}
{(ab)^n(2n(n+2)ab-(a^2+b^2)(n+1)^2)+a^{2(n+1)} + b^{2(n+1)}}, \\
\\
C=\frac{(ab)^{n-1}( (n+1)(a^{n+3}-ab^{n+2} -3ba^{n+2}) +
((n-1)b^{n+3} +2(n+2)b^2a^{n+1}))}
{(ab)^n(2n(n+2)ab-(a^2+b^2)(n+1)^2)+a^{2(n+1)} +
b^{2(n+1)}}, \\
\\
D=\frac{(ab)^n( b^{n+1}(n-b(n-2))-2a^nb^2(n+1)-na^{n+1}(a-3b))
}{(ab)^n(2n(n+2)ab-(a^2+b^2)(n+1)^2)+a^{2(n+1)} + b^{2(n+1)}}
\end{array}
\end{equation}
(cf. (3.13) in \cite{cal82:ext_kah_met}). Substituting \eqref{soln}
into \eqref{text} and \eqref{symp} gives the (second derivative of
the) $t$-potential $F_{\widehat{\bbC\bbP}^n}$ and the function
$h_{\widehat{\bbC\bbP}^n}$ of the extremal K\"ahler metrics,
respectively. We conclude that
\begin{coro}
There exists a family of $\U(n)$-invariant, extremal K\"ahler
metrics on $\widehat{\bbC\bbP}^n$ and these are determined by the
symplectic potential \eqref{symp} on $\Delta_{\widehat{\bbC\bbP}^n}$
where $h_{\widehat{\bbC\bbP}^n}(t) $ is a smooth function on
$[a,b]\in (0,\infty)$ given by \eqref{symp2} (as determined by
\eqref{soln}), and the real numbers $a,b$ parameterize the
cohomology class of these extremal K\"ahler metrics. Furthermore,
the scalar curvature of these extremal K\"ahler metrics is
$S(g_{\widehat{\bbC\bbP}^n})=At+B$ with $A,B$ as given by
\eqref{soln}.
\end{coro}
\begin{example}
{\rm In \S 6 of \cite{abr98:kah_geom_tor_var_ext_met} Abreu
derived a formula for the symplectic potential of Calabi's
extremal metric on $\widehat{\bbC\bbP}^2$ by Legendre transforming
the K\"ahler potential (as derived by Calabi) for this extremal
K\"ahler metric into symplectic coordinates. By his conventions
$0<a<1$, $b=1$ and $n=2$ i.e. $a$ is the amount by which
$\widehat{\bbC\bbP}^2$ is blown-up. Setting these values in
\eqref{soln} gives
\[
\begin{array}{c}
A = \frac{-24a}{a^3+3a^2-3a-1}, \;B =
\frac{6(3a^2-1)}{a^3+3a^2-3a-1}, \;C =
\frac{(3a^2-1)a}{a^3+3a^2-3a-1}, \;D =\frac{-2a^3}{a^3+3a^2-3a-1}.
\end{array}
\]
Substituting these into \eqref{symp2} gives
\[
h''_{\widehat{\bbC\bbP}^2}(t)=\frac{2a(1-a)}{2at^2+t-a^2t+2at+2a^2}
-\frac{1}{t}
\]
which is exactly Abreu's equation (24) in
\cite{abr98:kah_geom_tor_var_ext_met}.}
\end{example}

\bibliography{mybib}

\end{document}